\documentclass[11pt,b5paper,twoside,headrule]{amsart}
\usepackage{amsmath}
\usepackage{amsfonts}
\usepackage{amssymb}
\usepackage{euscript}
\usepackage{graphicx}
\footskip=10mm\parskip 0.2cm \addtocounter{page}{0}
\newtheorem{theorem}{\sc{Theorem}}

\newcounter{const}
\renewcommand{\theconst}
{D_{\arabic{const}}}
\newcommand{\newconst}
{\addtocounter{const}{1} \theconst}
\newcounter{con}
\renewcommand{\thecon}
{d_{\arabic{con}}}
\newcommand{\newcon}
{\addtocounter{con}{1} \thecon}
\newcommand{\e}{\varepsilon}

\newcommand{\T}{\mathbb{T}}

\newcommand{\F}{\EuScript{F}}
\begin{document}
\font\tenhtxt=eufm10 scaled \magstep0 \font\tenBbb=msbm10 scaled
\magstep0 \font\tenrm=cmr10 scaled \magstep0 \font\tenbf=cmb10
scaled \magstep0
\def\evenhead{{\protect\centerline{\textsl{\large{Ya.G. Sinai \quad and \quad M.D. Arnold }}}\hfill}}

\def\oddhead{{\protect\centerline{\textsl{\large{Global Existence and Uniqueness Theorem ...}}}\hfill}}

\pagestyle{myheadings} \markboth{\evenhead}{\oddhead}
\thispagestyle{empty}
\noindent{{\small\rm Pure and Applied Mathematics Quarterly\\ Volume 4, Number 1\\
(\textit{Special Issue: In honor  of \\Gregory Margulis, Part 2 of 2})\\
1---9, 2008} \vspace*{1.5cm} \normalsize

\begin{center}
{\bf{\Large Global Existence and Uniqueness Theorem for $3D$ --
Navier-Stokes System on $\T^3$ for Small Initial Conditions in the
Spaces $\Phi(\alpha).$}}
\end{center}

\begin{center}

\large{Ya.G. Sinai \quad and \quad M.D. Arnold}
\end{center}

\begin{center}
{\it Dedicated to G.A. Margulis on the occasion of his sixtieth
birthday.}
\end{center}

\footnotetext{Received October 11, 2005.}
\bigskip

\begin{center}
\begin{minipage}{5in}

\noindent \textbf{Abstract:}We consider Cauchy problem for
three-dimensional Navier-Stokes system with periodic boundary
conditions with initial data from the space of pseudo-measures
$\Phi(\alpha)$. We provide global existence and uniqueness of the
solution for sufficiently small initial data.
\end{minipage}
\end{center}

\section{Introduction}

Three-dimensional Navier-Stokes system with periodic boundary
conditions after Fourier transform can be written in the form:
\begin{equation}
\label{eq: NSS_velo} v(t,k)\! =\! \exp{-t|k|^2}v_0(k) + 2\pi i
\!\!\!\int\limits^t_0 \exp\{-(t - s)|k|^2\}\!\!\!\sum\limits_{l\in
\mathbb{Z}^3}\langle k, v(s,k-l)\rangle P_k v(s,l)ds
\end{equation}
Here $k\in \mathbb{Z}^3$, $t\in \mathbb{R}_+$, $v(t,k) \in
\mathbb{C}^3$, $v(t,k)\bot k$ for any $k\ne 0$ and $v(t,0) = 0$
for all $t> 0$. $v_0(k)$ is the initial condition and $P_k$
denotes Leray projector to the subspace orthogonal to $k$ and has
the form $P_k=\textrm{Id}-\dfrac{\langle k,
\cdot\rangle}{|k|^2}k$. Also \eqref{eq: NSS_velo} assumes that the
viscosity $\nu = 1$ and that the external forcing is absent.

T. Kato in \cite{K} proved the local existence theorem for the
$3D$–-Navier-Stokes system on $\mathbb{R}^3$ and global existence
and uniqueness theorem in the space $L^{\frac
32}(\mathbb{R}^3)\cap L^1(\mathbb{R}^3)$ for small initial
conditions.

In this paper we consider Cauchy problem for the system \eqref{eq:
NSS_velo} with initial data from the space $\Phi(\alpha)$ which is
analogous to the subspace $\Phi(\alpha,\alpha)$ introduced in
\cite{S1}, \cite{S2} and consists of functions of the form

\[
\Phi(\alpha)=\left\{f(k) = \frac{c(k)}{|k|^\alpha}, k\ne 0\mid
\sup\limits_{k}|c(k)|<\infty\right\}, \qquad
\|f(k)\|_{\alpha}=\sup\limits_{k\in \mathbb{Z}^3}|k|^\alpha |f(k)|
\]
We assume $\alpha>2$ and shall write $\alpha=2+\varepsilon$. V.
Kaloshin and Yu. Sannikov announced the global existence theorem in
the spaces $\Phi(\alpha)$, $\alpha\geqslant 2$ for small initial
data (see \cite{KS}). In this paper we give a detailed proof of this
result which shows also the character of decay of solutions in this
case. It is worthwhile to mention that according to our point of
view a similar result is not valid in the continuous case of $k\in
\mathbb{R}^3$.

Second author acknowledges financial support from NSF Grant DMS
0600996

\section{Main result}
The purpose of this paper is to prove the following theorem.
\begin{theorem}
\label{th: Main theorem}  Let $0 < 3\varepsilon< 1$ and
$\|v_0\|_\alpha\leqslant \delta$ where $v_0 =
\dfrac{c_0(k)}{|k|^\alpha}$ is the initial condition and
$\delta=\delta(\alpha)$ is sufficiently small. Then the equation
\eqref{eq: NSS_velo} has a global solution $v(t,k) = \dfrac{c(t,
k)}{|k|^\alpha}$ such that $c(t, k)$ is a continuous mapping of
$[0,\infty)$ into $L^\infty(\mathbb{Z}^3 \setminus \{0\})$, $t>0$.
\end{theorem}

The proof of the Theorem \ref{th: Main theorem} goes by induction.
Put $H^{(0)}_0(k) = \dfrac{c_0(k)}{|k|^\alpha}$, $k\ne 0$,
$H^{(1)}_0(k)=G_0(k)=0$ and assume that for some integer $m$ we
constructed the solution $v(t, k)$, $0 \leqslant t \leqslant m$,
such that
\begin{equation}
\label{eq: velo_Inductive} v(m, k) = H^{(0)}_m (k) + H^{(1)}_m (k)
+ G_m(k) \end{equation}

\noindent where
\begin{equation*}
\label{eq: H_0} H^{(0)}_m (k) =
\frac{\exp\{-m|k|^2\}c_0(k)}{|k|^\alpha},\end{equation*}

\begin{equation*}
\label{eq: H_1} H^{(1)}_m (k) = \sum\limits_{j=1}^m\exp\{-(m -
j)|k|^2\}h^{(1)}_j(k)
\end{equation*}
and
\begin{equation*}
\label{eq: G} G_m(k) =\sum\limits_{j=1}^m \exp\{-(m -
j)|k|^2\}g_j(k)
\end{equation*}

\noindent Suppose that for all $j \leqslant m$ functions
$h^{(1)}_j (k)$ satisfy the inequalities:
\begin{equation}
\label{eq: estimate h_1j} |h^{(1)}_j (k)| \leqslant
\frac{\newconst \delta^2\exp\left\{-\frac j 2
|k|^2\right\}}{|k|^{2\varepsilon}} ,
\end{equation}

\noindent while the  functions $g_j(k)$ satisfy the inequalities:

\begin{equation}\label{eq: estimate g_j}
|g_j(k)| \leqslant \frac{\newconst\delta^2\exp\{-\newcon
|k|\sqrt{m}\}}{|k|^\beta}
\end{equation}

\noindent Here  $\beta > 3$ is a constant and $d$, $D$ with
indices denote various absolute constants which appear during the
proof, but their exact values play no role in the arguments.

Consider $0 \leqslant t \leqslant 1$ and write down the solution
of \eqref{eq: NSS_velo} in the form:

\begin{equation}
\label{eq: inductive assumption} \begin{split}v(t + m, k) =
\frac{\exp\{-(m + t)|k|^2\}c_0(k)}{|k|^\alpha}\!\! +
\!\!\sum\limits_{j=1}^m \frac{\exp\{-(m - j +
t)|k|^2\}h^{(1)}_j(k)}{|k|^{2\e}} +\\+ \frac{h^{(1)}_{m+1}(t,
k)}{|k|^{2\e}} + \sum\limits_{j=1}^m \exp\{-(m - j + t)|k|^2\}g_j(k)
+ g_{m+1}(t, k)\end{split}\end{equation}

\noindent We show that the inequalities \eqref{eq: estimate h_1j},
\eqref{eq: estimate g_j} holds for $h^{(1)}_{m+1}(1,k)$ and
$g_{m+1}(1,k)$ respectively.

\section{Proof of the main result.}
Denote

\begin{equation}
\label{eq: H_0(t)}H^{(0)}_{m+1}(t, k)= \dfrac{\exp\{-(m +
t)|k|^2\}c_0(k)}{|k|^\alpha},\end{equation}

\begin{equation}
\label{eq: H_1_m+1} H^{(1)}_{m+1}(t, k)=\sum\limits_{j=1}^m
\frac{\exp\{-(m-j+t)|k|^2\}h^{(1)}_j(t, k)}{|k|^{2\e}}+
\frac{h^{(1)}_{m+1}(t, k)}{|k|^{2\e}},
\end{equation}

\begin{equation}
\label{eq: G_sum_m+1} G_{m+1}(t, k) = \sum \limits_{j=1}^m
\exp\{-(m - j + t)|k|^2\}g_j(t, k)
\end{equation}

\noindent and
\[
\left(H' \circledast H''\right)(t,k) = i \int\limits_0^t \exp\{-(t
- s)|k|^2\}\sum\limits_{k\in\mathbb{Z}^3\setminus
\{0\}\atop{k-l\ne 0}} \frac{\langle k,H'(s, k - l)\rangle P_k
H''(s, l)}{|k - l|^\alpha |l|^\alpha}\]

\noindent If we substitute \eqref{eq: inductive assumption} into
\eqref{eq: NSS_velo} we can write the expression for
$h^{(1)}_{m+1}(t, k)$:
\begin{equation}\label{eq: h_1_m+1}
h^{(1)}_{m+1}(t, k) = |k|^{2\e}\left(H^{(0)}_{m+1}\circledast
H^{(0)}_{m+1}\right)(t,k)
\end{equation} and the expression for
$g_{m+1}(t, k)$:

\begin{equation}\label{eq: g_m+1}
g_{m+1}(t, k) =\sum\limits_{j_1=1}^8 I^{(1,j_1)}_{m+1}(t, k)
+\sum\limits_{j_2=1}^3 I^{(2,j_2)}_{m+1} (t, k) + I^{(3)}_{m+1}(t,
k)
\end{equation}
where
\[\begin{split}
I^{(1,j_1)}_{m+1}(t, k) &= \left(H'\circledast H''\right)(t,k),\\
I^{(2,j_2)}_{m+1}(t, k) &= \left(H'\circledast
g_{m+1}\right)(t,k)+\left(g_{m+1}\circledast H'\right)(t,k)
\end{split}
\]
and $H'$, $H''$ are either $H^{(0)}_{m+1}(t, k)$, or
$H^{(1)}_{m+1}(t, k)$ or $G_{m+1}(t, k)$ except the case
$H'=H''=H^{(0)}_{m+1}$ which corresponds to the $h^{(1)}_{m+1}$
according to \eqref{eq: h_1_m+1}. Therefore $j_1$ changes from $1$
to $8$ and $j_2$ changes from $1$ to $3$. Also
\[
I^{(3)}_{m+1}(t, k) = g_{m+1}\circledast g_{m+1}.
\]

We see that $I^{(1)}_{m+1}$ does not depend on $g_{m+1}(t, k)$,
$I^{(2)}_{m+1}$ is a linear function of $g_{m+1}(t, k)$ and
$I^{(3)}_{m+1}$ is a quadratic function of $g_{m+1}(t, k)$.
Therefore \eqref{eq: g_m+1} is a typical equation which can be
solved by iterations if the coefficients are small enough. Below
we provide necessary estimates and later we come back to the
analysis of (8).

\subsection{First estimates.}

\noindent Here we show that all functions $h_{m+1}^{(1)}$ behaves
like gaussian functions of $t|k|$ and then provide necessary
estimates for coefficients in \eqref{eq: g_m+1}.

 As in \cite{S1}, \cite{S2}, we use the identity:
\begin{equation}
\label{eq: identity} a_1|k - l|^2 + a_2|l|^2 = \frac{a_1a_2}{a_1 +
a_2} |k|^2 + (a_1 + a_2)\left|l-\frac{a_1}{a_1+a_2}k\right|^2
\end{equation}

{\underline{\bf An estimate of $H_{m+1}^{(1)}$.}} At first we
estimate $h^{(1)}_{m+1}(t,k)$. From \eqref{eq: h_1_m+1} it follows
that

\[
\begin{split}
h_{m+1}^{(1)}(t,k)=&(H^{(0)}_{m+1}\circledast H^{(0)}_{m+1})(t, k)
= 2\pi i \int\limits^t_0 \exp\{-(t - s)|k|^2\}\cdot \\&\cdot
\sum\limits_{l\in \mathbb{Z}^3\setminus\{0\} \atop{k-l\ne 0}}
\frac{\langle k, c_0(k - l)\rangle P_k
c_0(l)}{|k-l|^\alpha|l|^\alpha} \exp\{-(m + s)|k-l|^2-(m + s)|l|^2
\}ds \end{split}
\]

\noindent Using \eqref{eq: identity} we can write
\[
\begin{split}
|h^{(1)}_{m+1}(t, k)| \leqslant &\exp\left\{-
\frac{m|k|^2}{2}\right\} \int\limits^t_0 \exp\{-(t - \frac
s2)|k|^2\} \cdot
\\
&\cdot \sum\limits_{l\in \mathbb{Z}^3\setminus\{0\}\atop{k-l\ne
0}} \frac{\langle k, c_0(k-l)\rangle P_kc_0(l)}{|k - l|^\alpha
|l|^\alpha} \exp\{-2m|l - \frac 12 k|^2\}ds \leqslant
\\
&\leqslant \delta^2 \exp\left\{- \frac{m|k|^2}{2}\right\}
\frac{\exp\{-\frac{t}{2} |k|^2\} - \exp\{-t|k|^2\}}{|k|^2}
\cdot\\&\cdot\sum\limits_{l\in\mathbb{Z}^3\setminus\{0\}\atop{k-l\ne0}}
\frac{\exp\{-2m|l - \frac 12 |k|^2\}}{|k -
l|^\alpha|l|^\alpha}\leqslant\\ &\leqslant
\frac{\newconst\delta^2}{|k|^{2\e}} \exp\left\{-\frac{(m +
t)|k|^2}{2}\right\}\frac{1 - \exp\{-\frac{t}{2}|k|^2 \}}{|k|^2}
\end{split}
\] Substituting this inequality
to \eqref{eq: H_1_m+1} we conclude

\begin{equation}
\label {eq: estimate H_1}
\begin{split}
\left|H^{(1)}_{m+1}(t, k)\right| \leqslant \frac{\theconst
\delta^2}{|k|^{2\e}} \frac{(1 -
\exp\{-\frac{t}{2}|k|^2\})}{|k|^2}\sum\limits_{j=1}^{m+1}
\exp\left\{ -(m + 1 - \frac j 2 )|k|^2\right\}\leqslant\\
\leqslant\frac{
\newconst \delta^2}{|k|^{2\e}}\frac{(1 -
\exp\{-\frac{t}{2}|k|^2\})}{|k|^2} \exp\left\{-\frac{(m +
1)|k|^2}{2}\right\}
\end{split}\end{equation}

{\underline{\bf Estimates for $H^{(j_1)}_{m+1}\circledast
H^{(j_2)}_{m+1}$.}} We present detailed estimate only for
$H^{(0)}_{m+1}\circledast H^{(1)}_{m+1}$ since all other terms can
be estimated in the same manner. From \eqref{eq: identity} and
\eqref{eq: estimate H_1} we have
\[\begin{split}
\left|(H^{(0)}_{m+1}\circledast H^{(1)}_{m+1})(t,k)\right|\leqslant
|k|\delta^3\int\limits^t_0 \exp\{-(t - s)|k|^2\}\cdot\\\cdot
\sum\limits_{l\in \mathbb{Z}^3\setminus\{0\}\atop{k-l\ne 0}}
\frac{\exp\{-(m+ s)|l|^2 - \frac{m +  s}{2} |k - l|^2\}}{|k -
l|^{2\e}|l|^\alpha} ds\leqslant
\end{split}\]
\[
\leqslant
|k|\delta^3\exp\left\{-\frac{m+1}{3}|k|^2\right\}\int\limits^t_0
\exp\{-(t - s)|k|^2\}\!\!\! \sum\limits_{l\in
\mathbb{Z}^3\setminus\{0\}\atop{k-l\ne 0}} \frac{\exp\{-\frac 32(m+
s)|l- \frac{m+s}{3}k|^2\}}{|k - l|^{2\e}|l|^\alpha} ds
\]
Since $\alpha+2\e> 2$ the last sum is not more than some constant
$\newconst$. We get

\begin{equation}\label{eq: estimate H0H1}\left|(H^{(0)}_{m+1}\circledast H^{(1)}_{m+1})(t,k)\right|
\leqslant |k|\theconst
\exp\left\{-\frac{m+1}{3}|k|^2\right\}\frac{1-\exp\{-t|k|^2\}}{|k|^2}\end{equation}

\noindent Similarly, for $(H^{(1)}_{m+1}\circledast
H^{(1)}_{m+1})(t,k)$ we can write

\begin{equation}
\label{eq: H1H1} \left|(H^{(1)}_{m+1}\circledast
H^{(1)}_{m+1})(t,k)\right| \leqslant \newconst
\exp\left\{-\frac{m+1}{4}|k|^2\right\}\frac{1-\exp\{-t|k|^2\}}{|k|}
\end{equation}

\subsection{Spaces $\F_m(c)$.} Fix positive constant $\beta>0$
and introduce functional space $\F_m(c)$
\[\F_m(c)=\left\{f(k)\mid |f(k)|\leqslant \frac{\newconst}{|k|^\beta}\exp\{-c\sqrt{m}|k|\}, \, k\ne 0 \right\},
\:\|f\|_{m,c}=\inf \theconst\]

We show that functions $g_{m+1}(t,k)$ belong to the spaces $\F_m$
with uniform constant if only all coefficients in \eqref{eq:
g_m+1} are sufficiently small. First of all we show that
$H^{(j_1)}\circledast H^{(j_2)}$, $j_1+j_2> 0$ belongs to the
space $\F_m(\newcon)$ for some constant $\thecon$. It follows from
previous estimates, that all of these functions decay as a
Gaussian functions. For our purpose it is convenient to consider
them as functions from the space $\F_m(\thecon)$. Since $m|k|
\geqslant 1$ we can write

\[\exp\left\{-\frac{m|k|^2}{3}\right\}\leqslant \frac{\newconst}{|k|^\beta}\exp\left
\{-\frac{\sqrt{m}|k|}{\sqrt{3}}\right\}\] for some constant
$\theconst$. We see that $(H^{(0)}_{m+1}\circledast
H^{(1)}_{m+1})(t, k) \in \F_{m+1}(\frac{1}{\sqrt{3}})$ and
\begin{equation}\label{eq: embedding}\|H^{(0)}_{m+1}\circledast
H^{(1)}_{m+1}\|_{m+1,\frac{1}{\sqrt{3}}} \leqslant
\newconst\end{equation} for some constant $\theconst$, which does not
depend on $t$.

Assuming that $G_{m+1}\in \F_{m+1}(\thecon)$ we can write for
$G_{m+1}\circledast H^{(0)}_{m+1}$

\[|(G_{m+1}\circledast H^{(0)}_{m+1})(t, k)| \leqslant \|G_{m+1}\|_{m+1,\thecon}
\delta |k| \int\limits^t_0 \exp\{-(t - s)|k|^2\}
\cdot\]\[\cdot\sum\limits_{l\in
\mathbb{Z}^3\setminus\{0\}\atop{k-l\ne 0}} \frac{ \exp\{-\thecon
\sqrt{m}|k-l| - m|l|^2\}}{|l|^\alpha |k-l|^{\beta_1}}ds\] For the
last expression we get
\[\begin{split}\exp\{-\thecon\sqrt{m}|k -l| - m|l|^2\}&\leqslant
\exp\{-\thecon\sqrt{m}|k|\} \exp\{\thecon\sqrt{m}|l| - m|l|^2\}
\leqslant\\&\leqslant\newconst \exp\{-\thecon\sqrt{m}|k|\}
\exp\{|l - \frac{\thecon}{2\sqrt{m}}|^2\}\end{split}\]
 So for
$(G_{m+1}\circledast H^{(0)}_{m+1})(t, k)$ we obtain
\[|(G_{m+1}\circledast H^{(0)}_{m+1}(t, k))| \leqslant\newconst
\|G_{m+1}\|_{m+1,\thecon}\delta
 \frac{\exp\{-\thecon|k|\sqrt{m}\}}{|k|^{\beta_1}}\frac{1 -
 \exp\{-t|k|^2\}}{|k|^2}\]
 All other terms in $I^{(1)}_{m+1}(t, k)$ can be similarly
estimated.

Thus we embed the first term in the representation of $g_{m+1}(t,
k)$ given by \eqref{eq: g_m+1} into the space $\F_{m+1}(\thecon)$.
Now we provide the necessary estimates for the terms
$I_{m+1}^{(3)}$ and $I^{(2)}_{m+1}$.

{\underline{ \bf Estimate for $I^{(3)}_{m+1}$.}} We show, that for
given functions $f_1$, $f_2 \in \F_{m+1}(\thecon)$ $f_1\circledast
f_2$ also belongs to the space $\F_{m+1}(\thecon)$ and $\|f_1,
f_2\|_{m+1,\thecon} \leqslant
\newconst\|f_1\|_{m+1,\thecon}\|f_2\|{m+1,\thecon}$ for some constant $\theconst$.

Write down the estimate
\[|f_1\circledast f_2|\leqslant \|f_1\|_{m+1,\thecon}\|f_2\|_{m+1,\thecon} |k|
\int\limits^t_0 \exp\{-(t - s)|k|^2\}\cdot\]\[\cdot
\sum\limits_{l\in \mathbb{Z}^3\setminus\{0\}\atop{k-l\ne 0}}
\frac{ \exp\{-\thecon \sqrt{m+1}(|l| - |k -l|\}}{|l|^\beta
|k-l|^\beta}ds\leqslant
\]\[
\leqslant \frac{\newconst
\|f_1\|_{m+1,\thecon}\|f_2\|_{m+1,\thecon}}{|k|^{2\beta-3}}\exp\{-\thecon|k|\sqrt{m+1}\}\frac{1-\exp\{-t|k|^2\}}
{|k|}
\]
Since $\beta > 3$ and the last expression is not more than $1$, we
get

\begin{equation}
\label{eq: contraction}|f_1\circledast f_2|\leqslant
\frac{\newconst
\|f_1\|_{m+1,\thecon}\|f_2\|_{m+1,\thecon}}{|k|^\beta}
\exp\{-\thecon\sqrt{m + 1}|k|\}
\end{equation} In particular, for $g_{m+1}(t, k) \in
\F_{m+1}(\thecon)$ it follows that
\begin{equation}\label{eq: estimate I3}\|I^{(3)}_{m+1}(t, k)\|_{m+1,\thecon}\leqslant
\theconst \|g_{m+1}(t, k)\|^2_{m+1,\thecon}\end{equation}

\paragraph{Estimates for $I^{(2)}_{m+1}$.}
Here we produce the upper bound  for $\|I^{(2)}(t,
k)\|_{m+1,\theconst} = \sum\limits_{j_2=1}^3 I^{(2,j_2)}_{m+1} (t,
k)$ assuming, that $g_{m+1}(t, k) \in \F_{m+1}(\thecon)$.
\[|(g_{m+1}(t, k)\circledast H^{(0)}_{m+1}(t, k))|\leqslant \|g_{m+1}\|_{m+1,\thecon}\delta |k|
\int\limits^t_0 \exp\{-(t - s)|k|^2\}\cdot\]\[\cdot
\sum\limits_{l\in \mathbb{Z}^3\setminus\{0\}\atop{k-l\ne 0}}
\frac{ \exp\{-\thecon \sqrt{m}|k -l|-m|l|^2\}}{|k-l|^\beta
|l|^\alpha}ds
\]

Again for the last expression holds
\[\exp\{-\thecon\sqrt{m}|k -l| - m|l|^2\}\leqslant \newconst \exp\{-\thecon\sqrt{m}|k|\}
\exp\{|l - \frac{\thecon}{2\sqrt{m}}|^2\}\] So for $(g_{m+1}(t,
k)\circledast H^{(0)}_{m+1}(t, k))$ we can write
\[|(g_{m+1}(t, k)\circledast H^{(0)}_{m+1}(t, k))| \leqslant\newconst \|g_{m+1}\|_{m+1,\thecon}\delta
\frac{1 - \exp\{-t|k|^2\}}{|k|}\frac{\exp\{-\thecon
\sqrt{m}|k|\}}{|k|^\beta}\] For the terms $g_{m+1}\circledast
G_{m+1}$ we can produce an appropriate estimate using \eqref{eq:
contraction}. All other terms in $I^{(2)}_{m+1}(t, k)$ can be
estimated in a similar way.

Collecting all present estimates we see that for some constant
$newcon$

\begin{equation}\label{eq: final bound}
\|g_{m+1}(t, k)\|_{m+1,\thecon}\leqslant  \newconst
\delta^2+\newconst\delta\|g_{m+1}\|_{m+1,\thecon} +
\newconst\|g_{m+1}\|^2_{m+1,\thecon}
\end{equation}
So for sufficiently small $\delta$ all coefficients in \eqref{eq:
final bound} are small and the equation \eqref{eq: g_m+1} can be
solved by iterations. The solution $g_{m+1}(t,k)$ belongs to
$\F_{m+1}(\thecon)$ and unique in this class of functions. Each
function $g_{m+1}(t, k)$ provides the unique solution $v(m + t,
k)$ of \eqref{eq: inductive assumption}. The Theorem \ref{th: Main
theorem} is proven.

\bigskip

\noindent Ya.G. Sinai\\Mathematics Department of Princeton
University, Princeton, NJ, USA\\
Landau Institute of Theoretical Physics, Moscow, Russia.\\

\noindent M.D. Arnold\\International Institute of the Earthquake
Prediction\\
Theory and Mathematical Geophysics, \\Russian Academy of Sciences,
Moscow, Russia.
\end{document}